\documentclass{article}
\usepackage{amssymb, amsmath}
\author{Shirley Bromberg \and Alberto Medina}
\date{February 22nd., 2001}
\title{G\'eom\'etrie des groupes oscillateurs et vari\'et\'es localement 
sym\'etriques}
\date{}
\def\dps{\displaystyle}
\def\r{{\mbox{\rm{I}\hspace{-0.02in}\rm{R}\hspace{0.005in}}}}
\def\C{\mbox{$\mathbb{C}$}}

\def\n{\mbox{$\mathbb{N}$}}
\def\qed{\hfill$\square$}
\def\g{{\cal{G}}}
\def\k{{\bf{k}}}

\def\G{{\bf{G}}}

\def\la{\lambda}
\def\dps{\displaystyle}
\def\be{\begin{equation}}
\def\ee{\end{equation}}
\newtheorem{theo}{Th\'eor\`eme}[section]
\newtheorem{pr}{Proposition}[section]
\newtheorem{lm}{Lemme}[section]
\newtheorem{df}{D\'efinition}[section]
\newtheorem{cor}{Corollaire}[section]

\begin{document}

\vspace{-0.8in}
\noindent
{\Large{\bf{G\'eom\'etrie des groupes oscillateurs et vari\'et\'es localement 
sym\'etriques
.}}}

\vspace{.2in}
\noindent
Shirley Bromberg

\noindent
E-mail address: stbs@xanum.uam.mx

\noindent
Departamento de Matem\'aticas, UAM-Iztapalapa, M\'exico, D.F., M\'exico.

\vspace{.1in}
\noindent
Alberto Medina

\noindent
E-mail address: medina@darboux.math.univ-montp2.fr

\noindent
D\'epartement des Math\'ematiques, 
Universit\'e de Montpellier II, France.

\vspace{.5in}
\noindent
\begin{minipage}[h]{11cm}
{\bf{Mot cl\'es:}}
Groupes Oscilateurs; Compl\'etude Geod\'esique; Structures 
(pseu\-do-riemanniennes) localement sym\'etriques; Vari\'et\'es localement 
affines.
\end{minipage}

\vspace{.3in}
\noindent
{\bf{MSC:}}
53C50,
53D25,
14R99, 
53C35.

\vspace{.15in}
\begin{minipage}[h]{10cm}
\footnotesize\baselineskip
\setlength{\baselineskip}
\centerline{{\bf{R\'esum\'e.}}}
Les groupes oscillateurs,$\G_\lambda,$ sont les seuls groupes de Lie simplement
connexes r\'esolubles non commutatifs qui admettent une m\'etrique de Lorentz
bi-invariante. Pour ces groupes nous:
d\'eterminons leurs groupes d'isom\'etries;
prouvons l'existence de structures affines invariantes \`a gauche et montrons 
que ces structures ne sont pas lorentziennes;
fournissons des conditions suffisantes pour que des m\'etriques 
pseudo-riemannienes invariantes \`a gauche soient compl\`etes.
Comme cons\'equence nous mettons en \'evidence des nouveaux exemples de vari\'et\'es
localement sym\'etriques parfois affines, parfois pseudo-riemanniennes (en
particulier lorentziennes) compactes compl\`etes ou incompl\`etes.
\end{minipage}

\vspace{.15in}
\begin{minipage}[h]{10cm}
\footnotesize\baselineskip
\setlength{\baselineskip}
\centerline{{\bf{Abstract.}}}
The Oscillator Groups,$\G_\lambda,$ are the only solvable, non commutative, simply 
connected Lie groups to admit a Lorentzian bi-invariant metric. 
For these groups,
we give sufficient conditions for a left-invariant pseudo-Riemannian
metric to be complete, we determine the group of isometries, we exhibit a 
left-invariant affine structure and prove that it is not
Lorentzian. As an application, we provide
new examples of compact pseudo-Riemannian (sometimes Lorentzian) 
locally symmetric, occasionally affine,
manifolds, complete or incomplete.
\end{minipage}

\section{Introduction-R\'esum\'e}

En 1985 Revoy et le deuxi\`eme auteur donnent dans \cite{kn:MR1} une 
collection tr\`es riche de nouveaux exemples de vari\'et\'es lorentziennes 
homog\`enes compactes. Ces vari\'et\'es sont des quotients de groupes de Lie,
appel\'es oscillateurs, qui sont
munis d'une m\'etrique de Lorentz bi-invariante,
par des sous-groupes de Lie co-compacts (r\'eseaux). Ces auteurs fournissent une
condition n\'ecessaire et suffisante pour qu'un
groupe oscillateur ait des r\'eseaux et classifient ceux-ci \`a isomorphisme 
pr\`es. Ces
r\'esultats sont la suite logique de leurs travaux sur
``la m\'ethode de la double extension ortogonale" qui permet de construire de 
fa\c{c}on
inductive les groupes de Lie appel\'es ``orthogonaux", c'est-\`a-dire qui 
admettent
une m\'etrique pseudo-riemanienne bi-invariante.

En 1988, D'Ambra en \cite{kn:da}, prouve que la composante connexe de 
l'identit\'e  ${\cal{I}}\,(M,g)_0$
du groupe des isom\'etries d'une vari\'et\'e lorentzienne analytique r\'eelle 
compacte simplement connexe est compacte.
Zeghib d'une part (\cite{kn:Ze}) et Adams et Stuck (\cite{kn:as1}) de l'autre, 
prouvent en 1997 que
les exemples de \cite{kn:MR0}, les vari\'et\'es consid\'er\'es par D'Ambra et les bons quotients
de ${\rm{SL}}(2,\r)$ muni de sa forme de Killing, sont ``essentiellement les 
seuls
exemples" (\cite{kn:Ze})
de vari\'et\'es lorentziennes homog\`enes de
volume fini dont le groupe des isom\'etries n'est pas compact.
Pour le faire ils utilisent des r\'esultats de
Zimmer (\cite{kn:zi}) et Gromov (\cite{kn:g}).

 Un des buts du pr\'esent travail est de fournir des nouveaux exemples de
vari\'et\'es localement sym\'etriques, parfois affines, parfois
pseudo-riemanniennes (en particulier lorentziennes) compactes compl\`etes ou
incompl\`etes. Pour le faire nous poursuivons l'\'etude de la g\'eom\'etrie des
groupes oscillateurs  entrepris dans \cite{kn:me}, \cite{kn:B-M} et \cite{kn:B-M2}.

Au paragraphe 3 nous donnons des conditions
suffisantes pour qu'une m\'etri\-que invariante \`a gauche sur un oscillateur 
soit compl\`ete (Th\'eor\`eme  \ref{theo:comp1}) puis nous
exhibons des structures de vari\'et\'e affine invariante \`a gauche sur ces groupes (Th\'eor\`eme 
\ref{theo:gaff}).

La proposition \ref{pr:naff}
affirme qu'aucune structure affine invariante \`a gauche sur $G_\lambda$ 
est la connexion de Levi-Civita d'une m\'etrique lorentzienne invariante.

De plus en dimension 4 aucune m\'etrique pseudo-riemannienne invariante \`a gauche, 
d'indice quelconque  est compatible 
avec une structure affine invariante (Proposition 
\ref{pr:os4}). La section s'ach\`eve en 
montrant l'existence de m\'etriques pseudo-riemanniennes invariantes \`a gauche,
d'indice quelconque non nul, sur 
les oscillateurs qui sont localement sym\'etriques et compl\`etes, donc 
sym\'etriques (Th\'eor\`eme \ref{theo:locsym}).
Ces m\'etriques sont en fait des g\'en\'eralisations de la structure 
orthogonale $\k_\lambda$ qui est unique \`a homothetie pr\`es 
(\cite{kn:MR2}).

Dans le paragraphe 4 nous determinons la composante connexe de l'unit\'e
${\rm{Isom}}\,(G,k)_0$ du groupe des isom\'etries d'un groupe de
Lie orthogonal $(G,k)$ avec $G$ simplement connexe (Th\'eor\`eme 
\ref{theo:acc}).
Nous specialisons ensuite
ce r\'esultat au cas d'un groupe oscillateur quelconque 
($G_\lambda,\k_\lambda).$

Puisque $(G_\lambda,\k_\lambda)$ est une vari\'et\'e lorentzienne localement 
sym\'etrique com\-pl\`ete, \`a tout sous-groupe $\Gamma$ de 
${\rm{Isom}}(G_\lambda,\k_\lambda)_0$ op\'erant librement et de
fa\c{c}on proprement discontinue sur $G_\lambda$ lui est associ\'ee 
une vari\'et\'e $G_\lambda/\Gamma=M_\lambda$
qui est naturellement munie d'une structure lorentzienne homog\`ene qui est
localement sym\'etrique et compl\`ete (car transitive) 
d'apr\`es le th\'eor\`eme classique \ref{theo:class}. Ceci g\'en\'eralise 
les exemples de \cite{kn:MR1}.

Le th\'eor\`eme \ref{theo:todo} munit les vari\'et\'es 
oscillatrices, 
respectivement de m\'etriques pseudo-riemanniennes localement
sym\'etriques compl\`etes (autres que celles d\'e\-dui\-tes de $\k_\lambda$),
des structures
affines compl\`etes ou non compl\`etes et des m\'etriques
pseudo-riemanniennes, d'indice quelconque non nul, non compl\`etes.

\section{Groupes Oscillateurs}

Pour $\lambda=(\lambda_1,\ldots,\lambda_n)$ dans $\r^n$ avec $0<\lambda_1\le\cdots\le\lambda_n,$
soit $G_\lambda$ le goupe de Lie de vari\'et\'e sous-jacente 
$\r^{2n+2}\equiv\r\times\r\times\C^n$ et produit 
$$
(t,s,z_1,\ldots,z_n)\cdot(t',s',z'_1,\ldots,z'_n)=
$$
$$
(t+t',s+s'+\frac{1}{2}\sum_{j=1}^n{\rm{Im}}\bar{z_j}\exp(i\, t\lambda_j)z'_j,
\ldots ,z_j+\exp(i\, t\lambda_j)z'_j,
\ldots).
$$

Les groupes $G_\lambda$ sont caract\'eris\'es par le fait d'\^etre les
seuls groupes de Lie simplement connexes, r\'esolubles et non commutatifs 
\`a admettre une m\'etrique de Lorentz
bi-invariante (\cite{kn:MR1}). De plus ils contiennent des
sous-groupes discrets co-compacts si et seulement si l'ensemble
$\{\lambda_1,\ldots,\lambda_n\}$
engendre un sous-groupe discret de $(\r,+)$ (\cite{kn:MR2}).

\begin{df}
Les groupes $G_\lambda$
sont appel\'es {\rm{Groupes Oscillateurs,}} leurs alg\`e\-bres de Lie,
${\rm{\g}}_\lambda,$ seront
dites {\rm{Alg\`ebres Oscillatrices}}. Dans le cas o\`u $\G_\lambda$
poss\`ede des sous-groupes discrets co-compacts $\Gamma,$ on appelle
{\rm{vari\'et\'e oscillatrice,}} une vari\'et\'e du type
$\G_\lambda/\Gamma.$
\end{df}

Les oscillateurs et le rev\^etement universel de ${\rm{SL}}(2,\r)$ 
(qui est le seul groupe r\'eel simple et simplement connexe
\`a admettre une m\'etrique de Lorentz bi-invariante)
ont la propri\'et\'e suivante (voir \cite{kn:Ze}):
Si $M$ est une vari\'et\'e lorentzienne homog\`ene de volume fini dont
le groupe des isom\'etries ${\cal{I}}(M)$ n'est pas compact
alors ou bien ${\cal{I}}(M)$ contient un rev\^etement fini de
${\rm{PSL}}(2,\r)$ ou bien ${\cal{I}}(M)$
contient un groupe oscillateur.

\vspace{.1in}

L'alg\`ebre de Lie $\g_\lambda$ du groupe de Lie $\G_\lambda$ 
est l'espace vectoriel
$$
\g_\lambda:=
{\rm{Vect}}\{e_{-1}, e_0,e_1,\ldots,e_n,\check{e}_1,\ldots,\check{e}_n\},
$$
muni du produit
$$
\begin{array}{ccc}
[e_{-1},e_j]=\lambda_j\check{e_j}&[e_j,\check{e_j}]=e_0&[e_{-1},\check{e_j}]=-\lambda_je_j
\end{array}
$$
(les autres crochets \'etant nuls ou d\'eduits par antisym\'etrisation).
Posons pour $x\in\g,$
$$
x=\sum_{j=-1}^n\, x_j\,e_j + \sum_{j=1}^n\, \check{x}_j\,\check{e}_j,
$$
La forme sym\'etrique non d\'eg\'en\'er\'ee, exprim\'ee dans
ces coordonn\'ees, 
$$
\k_\la(x,x):=2x_{-1}x_0+\sum_{j=1}^{n}\frac{1}{\la_j}\left(x_j^2+\check{x}_j^2\right)
$$
est orthogonale, c'est-\`a-dire, elle d\'efinit une m\'etrique
bi-invariante sur $\G_\lambda.$

Rappelons qu'il n'y a qu'une
seule structure orthogonale, \`a homoth\'etie
pr\`es, sur $\G_\lambda$ (voir ~\cite{kn:MR2}).

\section{M\'etriques pseudo-riemannienes et structures lo\-ca\-lement 
sy\-m\'e\-tri\-ques
in\-va\-rian\-tes sur les \newline groupes oscillateurs}

\'Etant donn\'e un groupe de Lie orthogonal $(\G,\k),$ il est clair que
les m\'etriques (pseudo-riemanniennes)  invariantes \`a gauche sur $\G$ sont en bijection avec
les isomorphismes $\k$-sym\'etriques de l'espace vectoriel $\g$ au moyen
de la formule:
\begin{equation}
\k_u(x,y):= \k(u(x),y)\label{eq:metr}.
\end{equation}

Soit $\k_u$ une m\'etrique invariante sur $\G$ et $\nabla^u$ la connexion
de Levi-Civita sur $\G$ associ\'ee. 

 Dans la suite, d'une part nous fournissons des conditions suffisantes pour
la compl\'etude g\'eod\'esique de la m\'etrique $\k_u$ 
puis nous mettons en \'evidence une structure
affine invariante \`a gauche sur $G_\lambda$ et nous entamons 
la question de la
platitude de $\k_u.$

Le fait que la courbe $t\mapsto\sigma(t)$ soit une g\'eod\'esique s'exprime
au niveau de $\g$ en disant que la courbe
$t\mapsto x(t):=(L_{\sigma(t)^{-1}})_{*,\sigma(t)}(\dot{\sigma}(t))$
est une solution de l'\'equation diff\'erentielle:
\be\label{eq:geo}
\dot{x}=-r(x, x)
\end{equation}
o\`u $r(x,y)$ est le produit de Levi-Civita sur $\g$ correspondant a
$\nabla^u.$ Rappelons la formule de Koszul suivante:
\be
\k_u(r(x, y),z)=\frac{1}{2}\left( \k_u([x,y],z)-\k_u([y,z],x)+\k_u([z,x],y)
\right).
\end{equation}
Que la connexion $\nabla^u$ soit g\'eod\'esiquement compl\`ete \'equivaut
\`a dire que le champ de vecteurs sur $\g$ d\'efini par (\ref{eq:geo})
est complet.

Si $\Phi_u$ est l'isomorphisme sym\'etrique de $\g$ dans $\g^*$ donn\'e
par
$$
\Phi_u(x)= \k(u(x), \cdot)
$$
la compl\'etude de (\ref{eq:geo}) revient \`a la compl\'etude du champ
dit d'Euler:
\begin{equation}\label{eq:Euler}
\dot{\xi}= - ad^*_{\Phi_u^{-1}\xi}\,\xi
\end{equation}
o\`u $ad^*$ d\'esigne la repr\'esentation co-adjointe de $\g.$
Le champ transport\'e par $\Phi_u^{-1}$ du champ d'Euler nous fournit
l'\'equation diff\'erentielle sur $\g$ suivante:
\begin{equation}\label{eq:geo2}
u(\dot{x})= [u(x),x]
\end{equation}
Le changement de variable $y=u(x)$ tranforme l'\'equation pr\'ec\'edente
en la {\it{paire de Lax}}:
\begin{equation}\label{eq:lax}
\dot{y}= [y,u^{-1}(y)].
\end{equation}
L'\'equation (\!\ref{eq:lax}) a deux int\'egrales premi\`eres:
$$
\begin{array}{ccc}
\k(y,y)&\mbox{et}&
\k(y,u^{-1}(y))
\end{array}
$$

\begin{theo}\label{theo:comp1}
Soit $u$ un isomorphisme $\k_\la$-sym\'etrique de l'espace $\g_\lambda.$
Si $u$ laisse stable ou bien le centre ${\rm{C}}(\g_\lambda)$ (et
\`a fortiori, son id\'eal d\'eriv\'e) ou bien une 
sous-alg\`ebre de Cartan de $\g_\la$
alors
la m\'etrique $\k_u$ sur ${\rm{G}}_\lambda$ est compl\`ete.
\end{theo}

\vspace{.1in}
\noindent
{\bf{Preuve.}} 
Si $u$ stabilise ${\rm{C}}(\g_\lambda)$, 
l'\'equation  (\ref{eq:Euler}) est \'equivalente \`a une \'equation
lin\'eaire; elle est donc compl\`ete.

\vspace{.1in}
Soit ${\cal{C}}$ une sous-alg\`ebre de Cartan. Puisque $(\g,\k_\la)$ est 
orthogonale, 
$\g={\cal{C}}\oplus {\cal{C}}^\perp$ (relat. \`a $\k_\la$) et comme $\g_\la$ 
est oscillatrice, avec $\k(E,eo)\ne 0,$ est un plan qui contient $E_0:=e_0$ 
et la restriction de $\k_\la$ est d\'efinie positive.  
${\cal{C}}={\rm{Ker\, adj}}_E$ est d'indice 1, o\`u $E$ est un \'el\'ement 
r\'egulier de $\g,$ c'est-\`a-dire
$\k(E,eo)\ne 0.$ Supposons que ${\cal{C}}$ soit laiss\'ee 
stable par $u$ alors l'espace euclidien ${\cal{C}}^\perp$ est aussi stable par
$u.$ Soit
$E_{1},\ldots,E_{2n},$ base $\k$ orthonorm\'ee  de ${\cal{C}}^{\perp}$ qui 
diagonalise $u\mid_{C^{\perp}}.$ Posons $u(E_i):=\mu_iE_i,$ avec
$\mu_1\leq\mu_2\leq\ldots\leq\mu_l< 0 <
\mu_{l+1}\leq\ldots\le\mu_{2n}.$  Soit $E_{-1}\in C$ tel que
$$
\k(E_{-1},E_{-1})=0,\;\;\; \k(E_{-1},e_0)=1.
$$
Les vecteurs $E_{-1},E_0=e_0,E_{1},\ldots,E_{2n}$ sont une base pour $\g.$ 
Posons pour $x\in\g,$ 
$x=\sum_{i=-1}^{2n}\, {\bar{x}}_iE_i.$ Comme $u$ est $\k$- sym\'etrique et 
pr\'eserve ${\cal{C}},$ 
$u(e_0)=a\,E_{-1} +\alpha\,e_0$ et $u(E_{-1})=\alpha\,E_{-1} +b\,e_0.$ 
Les int\'egrales premi\`eres sont donn\'ees par:
$$
\begin{array}{ccccccc}
E(x)&=&\k_u(x,x)&=&a{\bar{x}}_0^2+b{\bar{x}}_{-1}^2+2\alpha\,{\bar{x}}_{-1}{\bar{x}}_0 &+& \sum_{i=1}^{2n}\, \mu_i{\bar{x}}_i^2\\
A(x)&=&\k_u(u(x),x)&=& 
2(\alpha{\bar{x}}_{-1}+ a{\bar{x}}_0)(b{\bar{x}}_{-1}+\alpha {\bar{x}}_0)
&+&\sum_{i=1}^{2n}\, \mu_i^2{\bar{x}}_i^2
\end{array}
$$
De plus, le centre apporte aussi une int\'egrale premi\`ere:
$$
C(x)=\k(x,u(e_0))=a{\bar{x}}_0 +  \alpha{\bar{x}}_{-1}.
$$
Si $a=0,$ $u$ pr\'eserve le centre et 
la m\'etrique $\k_u$ est donc compl\`ete. Supposons $a\neq
0.$ Les expressions:
\begin{equation}\label{eq:int}
\frac{ab-\alpha^2}{a\mu_j}(\mu_j{\bar{x}}_{-1}- C(x))^2
+\sum_{i=1}^{2n}\, \mu_i(\mu_j-\mu_i){\bar{x}}_i^2
\end{equation}
($j=1,2,\ldots,2n$) sont des int\'egrales premi\`eres .
Une des int\'egrales premi\`eres qui correspondent \`a $j=l$ et $j=l+1$ est
d\'efinie n\'egative. Ceci implique que 
$\mu_j{\bar{x}}_{-1}- C(x)\; (j=l \, {\mbox{ou}}\, j=l+1)$  et $x_i$  
$(i=1,\ldots, 2n)$  
sont born\'ees. D'o\`u la compl\'etude.

\qed

\begin{theo}\label{theo:gaff}
 Le groupe oscillateur $\G_\lambda$ est muni d'une structure de vari\'et\'e
affine compl\`ete d\'efinie pour les champs invariants \`a gauche sur $\G_\lambda$ par
le produit $r$ sur $\g_\lambda$ donn\'e par
$$
\left\{
\begin{array}{lcl}
r(x,y)&=&(1/2)[x,y]\\
r(e_{-1},y)&=&[e_{-1},y]\\
r(x,e_{-1})&=&r(e_{-1},e_{-1})=0
\end{array}
\right.
$$
pour $x$ et $y$ dans l'id\'eal d\'eriv\'e de $\g_\lambda.$
\end{theo}

\noindent
{\bf{Preuve.}}
On constate directement que $r$ est un produit \`a associateur 
sym\'etrique \`a gauche compatible avec le crochet de $\g_\lambda.$ 
Ceci signifie que la connexion
invariante \`a gauche d\'efinie par $r$ est \`a courbure et torsion nulles.
Cette
connexion est compl\`ete car les multiplications droites pour $r$ sont 
nilpotentes. 

Il se pose alors la question de savoir si cette connexion de courbure et
torsion nulles, ou toute autre ayant ces propriet\'es, est compatible avec une
m\'etrique invariante \`a gauche.

Nous avons le r\'esultat suivant

\begin{pr}\label{pr:naff}
Aucune structure affine invariante \`a gauche sur $\G_\lambda$ est
compatible avec une m\'etrique lorentzienne invariante \`a gauche.
\end{pr}

\noindent
{\bf{Preuve.}} 
Supposons qu'il existe une m\'etrique lorentzienne, invariante \`a gauche, 
plate $<\, ,\,>$
sur $\G_\lambda.$ Il existe alors
un isomorphisme $u$ de l'espace $\g_\lambda$  $\k$-sym\'etrique
tel que $\k_u=<\, ,\,>.$ Alors $L_{e_0}=0,$  $e_0$ est vecteur propre de 
$u,$ et par cons\'equent $k_u$ isotrope.
Posons
$u(e_0)=\nu e_0.$   Comme $u$ laisse
stable ${\cal{Z}}(\g),$ il laisse stable $\g'$ et induit sur
$\g'/{\cal{Z}}(\g):=\tilde{\g}$ un isomorphisme $\tilde{u}.$ Les 
m\'etriques $\k$ et $\k_u$ induisent des m\'etriques 
$\tilde{\k}$ et $\tilde{\k}_u$ sur $\tilde{\g}$ 
d\'efinies positives. Soit ${\cal{R}}$ la courbure de la connexion 
associ\'ee \`a $\k_u.$ La condition $u(e_0)\in\r e_0$ implique, 
pour $x,y\in\g',$ ${\cal{R}}_{xy}=0$ et, pour $x,y\in\g,$ $L_xy\in\g'.$ 
Ainsi il nous suffit d'\'etudier ${\cal{R}}_{e_{-1}x}$ pour $x\in\g'.$

Dans ce but, il sera n\'ecessaire de faire quelques calculs. D'abord, 
il existe une base $\tilde{\k}$ orthonormale,
$\{\tilde{E}_i\},$  qui diagonalise  $\tilde{u}:$
$\tilde{u}(\tilde{E}_i)=\nu_i\tilde{E}_i.$ 

Soit $E_i\in\g'$ tel que sa projection sur $\tilde{\g}$ est $\tilde{E}_i.$ Alors
$\{e_0, E_1,\cdots, E_{2n}\}$ est base de $\g'$ et
$$
\k(E_i,E_j)=\delta_{ij}\quad\quad \k(e_0,E_i)=0\quad i,j\ge 1.
$$
\'Ecrivons
$$
u(E_i)=\nu_iE_i+\mu_ie_0,
$$
et posons pour $i,j\ge 1,$ $[E_i,E_j]=\rho_{ij}e_0.$
 Remarquons que $\rho_{ij}=-\rho_{ji}.$

Nous nous proposons de calculer
$$
<{\cal{R}}_{e_{-1}E_i}E_i,e_{-1}>.
$$
Il est facile de v\'erifier que
$[e_{-1},E_i]=\sum_j\rho_{ij}E_j$ et que 
$$<L_{E_j}e_{-1},E_i>=(\rho_{ij}/2)(\nu_i-\nu_j+\nu).
$$
Ainsi
$$
<{\cal{R}}_{e_{-1}E_i}E_i,e_{-1}>=\sum_j\rho^2_{ij}(\nu_j-\nu_i)
-\sum_j\frac{\rho^2_{ij}}{4\nu_j}(\nu_j-\nu_i+\nu)^2
$$
et
$$
\sum_i<{\cal{R}}_{e_{-1}E_i}E_i,e_{-1}>=
-\sum_{i,j}\frac{\rho^2_{ij}}{4\nu_j}(\nu_j-\nu_i+\nu)^2.
$$
La m\'etrique \'etant suppos\'ee lorentzienne et plate, $\nu_j>0$ et
$$
0=
\sum_{i,j}\frac{\rho^2_{ij}}{4\nu_j}(\nu_j-\nu_i+\nu)^2,
$$
implique pour  $i,j\ge 1,$
$$
0=
\rho^2_{ij}(\nu_j-\nu_i+\nu)^2.
$$
Or ils existent $i,j$ tels que $\rho_{ij}\ne 0,$ et par cons\'equent
$$
\nu_j-\nu_i+\nu=0\quad\quad
{\mbox{et}}\quad\quad
\nu_i-\nu_j+\nu=0
$$
d'o\`u $\nu=0,$ ce qui est impossible.
\qed

On peut se demander si aucune structure affine invariante \`a gauche sur 
$\G_\lambda$ est compatible avec une m\'etrique invariante \`a gauche. 
\`A ce propos nous avons 
la r\'eponse partielle suivante:

\begin{pr}[\cite{kn:B-M2}]\label{pr:os4}
Le groupe oscillateur de dimension 4 n'admet pas de m\'etrique
invariante plate.  
\end{pr}

Mais les groupes oscillateurs admettent plusieures m\'etriques 
invariantes \`a gauche sym\'etriques. En effet,


\begin{theo}\label{theo:locsym}
Les m\'etriques invariantes \`a gauche de $\G_\la,$ d\'efinies 
par les isomorphismes $\k$-sym\'etriques $u$ de l'espace
$\g_\la$ donn\'es par $u(e_0)=e_0,$  
$u(e_i)=\eta_ie_i,$ 
$u(\check{e}_i)=\check{\eta}_i\check{e}_i$ et qui v\'erifient
$$
\eta_i+\check{\eta}_i=1\eqno(a)
$$ ou bien $$\eta_i=\check{\eta}_i\eqno(b)$$ 
pour chaque $i$ compris entre 1 et $n$
sont sym\'etriques car elles sont localement sym\'e\-triques et 
compl\`etes pour tout $\la.$
\end{theo}
\vspace{.1in}
\noindent
{\bf{Preuve.}} Puisque $u$ stabilise le centre de $\g_\la,$
la m\'etrique $\k_u$ compl\`ete 
(Th\'eor\`eme \ref{theo:comp1}).

Si l'on pose $r(x,y)=L_xy$ et ${\cal{R}}(x,y)=L_{[x,y]}-[L_x,L_y]$ (
le tenseur de courbure),
que $k_u$ soit localement sym\'etrique s'exprime par l'identit\'e
\begin{equation}\label{eq:loc}
[L_z,{\cal{R}}(x,y)]={\cal{R}}(L_zx,y)+{\cal{R}}(x,L_zy),
\end{equation}
pour tous $x,y,z\in\g.$ 

Or, $L_x=(1/2)\left({\rm{ad}}_x-u^{-1}{\rm{ad}}_{u(x)}+
u^{-1}{\rm{ad}}_{x}u\right).$
Ainsi $u(e_0)=e_0$ implique
$L_{e_0}=0$ et $L_xy\in\g',\,\forall x,y\in\g;$
$L_xy\in\r e_0, [L_x,L_y]=0,{\cal{R}}(x,y)=0, 
\forall x,y\in\g';$
$L_{e_i}e_j=L_{\check{e}_i}\check{e}_j=0\quad\forall i,j;$
$L_{\check{e}_i}e_j=L_{e_i}\check{e}_j=0\quad\forall i\ne j.$

Tandis que les autres directions propres de $u$ donnent
$u(e_{-1})=e_{-1}+\rho e_0,$ $L_{e_{-1}}e_{-1}=0,$ 
$L_{\check{e}_i}e_i=-(1/2)(\eta_i-\check{\eta}_i+1)e_0,$ 
$L_{e_i}\check{e}_i=-(1/2)(\eta_i-\check{\eta}_i-1)e_0,$
et
$$
\begin{array}{cc}
L_{e_{-1}}e_i=\dps\frac{\la_i}{2\check{\eta}_i}
(\eta_i+\check{\eta}_i-1)\,\check{e}_i&
L_{e_i}e_{-1}=\dps\frac{\la_i}{2\check{\eta}_i}
(\eta_i-\check{\eta}_i-1)\,\check{e}_i\\
&\\
L_{e_{-1}}\check{e}_i=\dps\frac{-\la_i}{2\eta_i}
(\eta_i+\check{\eta}_i-1)\,e_i&
L_{\check{e}_i}e_{-1}=\dps\frac{\la_i}{2\eta_i}
(\eta_i-\check{\eta}_i+1)\,e_i.
\end{array}
$$

Si $i$ se trouve dans la cas (a), nous avons
$$
L_{e_{-1}}e_i=L_{e_{-1}}\check{e}_i=0,\quad
L_{e_i}e_{-1}=[e_{-1},e_i],\quad
L_{\check{e}_i}e_{-1}=[e_{-1},\check{e}_i],
$$
pendant que les conditions (b)
fournissent
$$
\begin{array}{cc}
L_{e_{-1}}e_i=\phantom{-}\dps\frac{\la_i}{2\eta_i}
(2\eta_i-1)\,\check{e}_i&
L_{e_i}e_{-1}=\dps\frac{-\la_i}{2\eta_i}\check{e}_i\\
&\\
L_{e_{-1}}\check{e}_i=\dps\frac{-\la_i}{2\eta_i}
(2\eta_i-1)\,e_i&
L_{\check{e}_i}e_{-1}=\dps\frac{\la_i}{2\eta_i}e_i.
\end{array}
$$

Pour $x,y\in\g'_\la,$ l'identit\'e (\ref{eq:loc})
est trivialement v\'erifi\'ee.
Reste \`a constater que l'on a,
$$
[L_z,{\cal{R}}(e_{-1},y)]={\cal{R}}(L_ze_{-1},y)+{\cal{R}}(e_{-1},L_zy)\\
$$
pour $y\in\g'_\la.$ Or dans ce cas ${\cal{R}}(L_ze_{-1},y)=0$ et 
(\ref{eq:loc}) devient,
$$
[L_z,{\cal{R}}(e_{-1},y)]={\cal{R}}(e_{-1},L_zy).
$$
Pour $z$ dans $\g_\la,$ ${\cal{R}}(e_{-1},L_zy)=0$ et $[L_z,{\cal{R}}(e_{-1},y)]
=-L_y\circ L_{e_{-1}}\circ L_z-L_z\circ L_{e_{-1}}\circ L_y.$ Ce dernier
\'etant nul sur $\g'_\la,$ et
$$
<L_y\circ L_{e_{-1}}\circ L_z e_{-1},e_{-1}>=
-<L_z\circ L_{e_{-1}}\circ L_y e_{-1},e_{-1}>,
$$
par antisym\'etrie des multiplications \`a gauche, l'\'equation est 
v\'erifi\'ee dans ce cas.

Soit $z=e_{-1}.$ Pour $i$ dans le cas (a), et $\xi=e_i$ ou 
$\xi=\check{e}_i$
$$
{\cal{R}}(e_{-1},L_{e_{-1}}\xi)=0
$$
et
$$
[L_{e_{-1}},{\cal{R}}(e_{-1},\xi)]=[L_{e_{-1}},L_{[e_{-1},\xi]}]-
[L_{e_{-1}},[L_{e_{-1}},L_{\xi}]]
$$
qui s'annule sur $\g_\la'$ et 
$$
[L_{e_{-1}},L_{[e_{-1},\xi]}]e_{-1}-
[L_{e_{-1}},[L_{e_{-1}},L_{\xi}]]e_{-1}=0={\cal{R}}(e_{-1},L_{e_{-1}}\xi)e_{-1}.$$
Finalement, pour $i$ dans le cas (b), on v\'erifie que
$$[L_{e_{-1}},{\cal{R}}(e_{-1},e_i)]e_{-1}={\cal{R}}(e_{-1},L_{e_{-1}}e_i)e_{-1}=
\frac{\la_i^3}{8\eta_i^3}(2\eta_i-1)\check{e}_i.
$$
\qed

\section{Isom\'etries de $(\G_\la,\k_\la)$}

Soit ${\cal{I}}:={\rm{Isom}}\,(G,g)$ le groupe des isom\'etries du groupe
$G$ muni d'une m\'etrique
pseudo-riemannienne invariante \`a gauche $g$ et soit 
${\cal{I}}_\varepsilon:={\rm{Isom}}_\varepsilon(G,g)$ le
sous-groupe ferm\'e du groupe ${\cal{I}}$ des isom\'etries qui fixent
l'\'el\'ement neutre $\varepsilon$ de $G.$ Le groupe ${\cal{I}}$ 
op\`ere
naturellement et transitivement sur $G.$ Conna\^\i\-tre ${\cal{I}}$ comme
vari\'et\'e revient \`a conna\^\i tre ${\cal{I}}_\epsilon.$ Or, 
tout \'el\'ement de ${\cal{I}}_\varepsilon$ est d\'etermin\'e par sa
diff\'erentielle en $\varepsilon,$ qui est un
\'el\'ement de ${\cal{O}}(\g,g)$ qui pr\'eserve la courbure. L'ensemble de
ces \'el\'ements, que l'on d\'esignera par ${\cal{O}}_{{\cal{R}}}(\g,g),$
est un sous-groupe ferm\'e de ${\cal{O}}(\g,g).$ On a en fait:

\begin{lm}
Si $(\G,\k)$ est un groupe orthogonal et $\G$ est simplement connexe, alors 
l'application
$$
\begin{array}{ccccc}
\Phi &:& {\cal{I}}_\varepsilon &\to& {\cal{O}}_{{\cal{R}}}(\g,g)\\
&&\phi&\mapsto &\phi_{*,\varepsilon}
\end{array}
$$
est un isomorphisme de groupes de Lie.
\end{lm}

\vspace{.1in}
\noindent
{\bf{Preuve.}}
Puisque $\k$ est bi-invariante, $(\G,\k)$ est complet et sym\'etrique
($\sigma\mapsto\sigma^{-1}$ est la sym\'etrie
de centre $\varepsilon$). En outre, pour $u\in{\cal{O}}_{{\cal{R}}}(\g,\k),$ 
l'isom\'etrie locale  
${\cal{P}}_u:={\rm{Exp}}_\k\circ u\circ {\rm{Log}}_\k,$ dite polaire de $u,$
s'\'etend en une unique isom\'etrie globale de $(\G,\k)$ car $\G$ est 
simplement connexe (voir, par exemple, \cite{kn:ON}

\qed

\begin{theo}\label{theo:acc}
Soit $(\G,\k)$ un groupe orthogonal, $\G$ simplement connexe. Alors
$
\left({\rm{Isom}}\,(\G,\k)\right)_0$ s'identifie \`a
la vari\'et\'e produit $\G\times\left({\cal{O}}_{{\cal{R}}}(\g,\k)\right)_0
$
munie du produit 
$(\sigma,u)\cdot(\sigma',u')=(\sigma (u\cdot\sigma'),(\sigma'\cdot u)\circ u')$
o\`u les actions 
de $\left({\cal{O}}_{{\cal{R}}}(\g,\k)\right)_0$ sur $\G$ et de $\G$ sur 
$\left({\cal{O}}_{{\cal{R}}}(\g,\k)\right)_0$ sont
donn\'ees respectivement par
$$
\begin{array}{ccl}
u\cdot\sigma&=&\Phi^{-1}(u)(\sigma)\\
\sigma\cdot u&=& \left(L_{u\cdot\sigma}\right)^{-1}_{*,u\cdot\sigma}
\circ u\circ \left(L_\sigma\right)_{*,\varepsilon},
\end{array}
$$
$L_\sigma$ \'etant la multiplication \`a gauche par $\sigma$
dans le groupe $G.$
\end{theo}

\vspace{.1in}
\noindent
{\bf{Preuve.}}
La vari\'et\'e produit $\G\times\left({\cal{O}}_{{\cal{R}}}(\g,\k)\right)_0$
s'identifie \`a ${\rm{Isom}}\,(\G,\k)$ par 
$(\sigma,u)\mapsto L_{\sigma}\circ {\cal{P}}_u.$ Comme ${\cal{P}}_u\circ 
L_{\sigma}$ est une isom\'etrie, il existe 
$\sigma'$ et
$u'$ uniques tels que
\begin{equation}\label{eq:pol}
{\cal{P}}_u\circ L_{\sigma}=L_{\sigma'}\circ {\cal{P}}_{u'}
\end{equation}
avec $u\cdot\sigma:=\sigma'= {\cal{P}}_{u}(\sigma)$ et 
$\sigma\cdot u:=u'= 
\left(L^{-1}_{{\cal{P}}_u(\sigma)}\circ{\cal{P}}_u\circ 
L_{\sigma}\right)_{*,\varepsilon}.$

\qed

\vspace{.2in}
Rappelons le th\'eor\`eme suivant qui motive le fait de d\'eterminer dans 
la suite le 
groupe d'isom\'etries de $(\G_\la,\k_\la)$  et en
particulier les sous-groupes qui op\`erent convenablement sur $G$ de sorte
\`a obtenir des varietes pseudo-riemanniennes localement sym\'etriques:

\begin{theo}[voir \cite{kn:w}, p. 62]\label{theo:class}
Les vari\'et\'es pseudo-riemanniennes com\-pl\`e\-tes connexes localement
sym\'etriques sont les quotients $M/\Gamma$ o\`u $\Gamma$ est
un groupe d'isom\'etries op\'erant librement et de fa\c{c}on proprement
discontinue sur une va\-ri\'e\-t\'e $M$ pseudo-riemannienne sym\'etrique 
simplement connexe.
\end{theo}

\vspace{.2in}
D\'eterminons pour $(\G_\la,\k_\la)$ les actions du th\'eor\`eme 
\ref{theo:acc}.

\vspace{.1in}
Soient
$r_1, r_2, \ldots, r_s\in\n$ tels que
$$
\lambda_1=\cdots=\lambda_{r_1}<
\lambda_{r_1+1}=\cdots=\lambda_{r_1+r_2}<\cdots<
\lambda_{n-{r_s}+1}=\cdots=\lambda_n,
$$
et pour $i=1,\cdots,s$
posons
$$V_i:={\rm{Vect}}\{ e_r,\check{e}_r;
\la_r=\la_i\}
\;\;\;\;\;\;
V_i^\sim:=\r e_0\oplus V_i.
$$
La restriction de $\k_\la$ \`a chacun des
sous-espaces
$V_i^\sim$ est d\'eg\'en\'er\'ee;
elle induit sur $V_i^\sim/\r e_0(\equiv V_i)$ une forme quadratique
euclidienne,  not\'ee
$\k_i.$ Nous avons
\begin{pr}\label{pr:lm1}
Le groupe ${\cal{O}}_{{\cal{R}}}(\g_\la,\k_\la)$ s'identifie au groupe produit 
direct du groupe multiplicatif ${\mathbb{Z}}_2$
avec les groupes $V_i\rtimes{\cal{O}}(V_i,\k_i)$ des
d\'eplacements rigides des espaces euclidiens $(V_i,\k_i),$ $i=1,\ldots ,s.$ 
\end{pr}

\vspace{.1in}
\noindent
{\bf{Preuve.}}
Soit $u\in{\cal{O}}_{{\cal{R}}}(\g_\la,\k_\la).$
Nous allons d\'emonter qu'il existe
$\rho\in {\mathbb{Z}}_2,$ et 
$(v_i, u_i)\in V_i\rtimes{\cal{O}}(V_i,\k_i),$ $i=i,\ldots , r,$
uniques tels que $u$ s'exprime comme:
\begin{eqnarray}
u(e_{-1})&=&\phantom{-}\rho e_{-1}+\alpha e_0+\sum_i v_i\label{eq:10}\\
u(\, e_{0}\,)&=&\phantom{-}\rho e_{0}\label{eq:11}
\end{eqnarray}
et pour $v\in V_i,$
$$
u(v)=-\rho\k_\la (u_i (v),v_i)e_0+u_i (v)
$$
o\`u  $\alpha=-\rho/2\sum_i\k_\la(v_i,v_i).$ 
Ceci implique, \'evidemment, le r\'esultat.

L'\'equation (\ref{eq:11}) \'equivaut \`a dire que
$e_0$ est vecteur propre de $u.$ Or, comme
${{\cal{R}}}_{xy}=(1/4){\rm{ad}}_{[x,y]}$ et $u$ pr\'eserve
le tenseur de courbure alors
$$
u([x,[y,z]])=[u(x),[u(y),u(z)]]\;\;\forall x,y,z\in\g_\lambda.
$$
Donc, pour tous $x,z\in\g_\lambda,$
$$
0=u([x,[e_0,z]])=[u(x),[u(e_0),u(z)]],
$$
d'o\`u
$u(e_0)\in{\rm{C}}^2(\g_\la)={\cal{Z}}(\g_\la)$ 
et
$u(e_0)=\rho e_0.$ D'autre part $u\in {\cal{O}}_{{\cal{R}}}(\g_\la,\k_\la)$
implique
$$
u(e_{-1})=(1/\rho) e_{-1}+\alpha e_0+\sum_{i=1}^r v_i,
$$
avec $v_i\in V_i.$ Comme $0=\k_\la(u(e_{-1}),u(e_{-1}))=(2\alpha/\rho)+\sum_i
\k_i(v_i,v_i),$ alors
$\alpha=-(\rho/2)\sum_i\k_i(v_i,v_i).$ 

L'isomorphisme 
$u$ induit un isomorphisme
$$
\tilde{u}\in{\cal{O}}(\g'/\r e_0, \tilde{\k}),
$$
o\`u $\tilde{k}$ est la forme quadratique d\'efinie positive
induite par $\k_\la.$
Posons
$$
u(e_i)=C_ie_0+\sum_{k=1}^n\left(A_{ki}e_k+B_{ki}\check{e}_k\right)
$$
Nous avons
$$
\phantom{!} [u(e_i),[u(e_{-1}),u(e_i)]]=
\dps\frac{1}{\rho}\sum_{k=1}^n\la_k\left(A_{ki}^2+B_{ki}^2\right)e_0
$$
Puisque $u$ pr\'eserve la courbure nous obtenons,
$$
\la_i\rho^2=
\dps\sum_{k=1}^n\la_k\left(A_{ki}^2+B_{ki}^2\right)
$$
pour tout $i=1,\ldots,n.$

Il nous faut d\'emontrer que $\rho^2=1.$
Supposons $\rho^2<1.$ Alors, 
$\rho^2\la_1^2<\la_k^2,$
pour tout $k=1,\ldots,n.$
Comme $u$ est orthogonale,
$$
\frac{1}{\la_i}=\sum_{k=1}^n\frac{1}{\la_k}\left(A_{ki}^2+B_{ki}^2\right).
$$
D'o\`u
\begin{equation}\label{eq:aquella}
\dps\sum_{k=1}^n\frac{\la_1^2\rho^2-\la_k^2}{\la_k}\left(A_{k1}^2+B_{k1}^2\right)
=0,
\end{equation}
et
$A_{j1}=B_{j1}=0.$
Ceci \'etant impossible, 
$\rho^2\ge 1.$
Un raisonement en tout analogue au pr\'ec\'edent d\'emontre que
$\rho^2\le 1.$ D'o\`u $\rho^2=1$ et l'\'equation (\ref{eq:aquella}) devient
$$
\dps\sum_{j>r_1}\frac{\la_i^2-\la_j^2}{\la_j}\left(A_{ji}^2+B_{ji}^2\right)=0
$$
pour tout $1\le i\le r_1.$ Donc
$$
A_{ji}=B_{ji}=0
$$
pour tout $j>r_1$ et tout $1\le i\le r_1,$ et il est de m\^eme pour
$\check{A}_{ji},\,\check{B}_{ji}.$ Par cons\'equent
$u$ laisse stable $V_1$ et, par recurrence, $u$ laisse stable
$V_i,$ $i=1,\ldots, s.$ 
Posons, pour $v\in V_i,$ $u(v)=: C(v)e_0+u_i(v).$ L'\'egalit\'e
$\k_\la(v,v)=\k_\la(u(v),u(v))=\k_\la(u_i(v),u_i(v)),$ montre que
$u_i\in {\cal{O}}(V_i,\k_i).$
Comme $u$ pr\'eserve $\k_\la,$ nous avons, pour $v\in V_i,$
$0=\k_\la (u(e_{-1}),u(v))=\rho C(v)+\k_\la(v,u_i(v)).$
D'o\`u $C(v)=-\rho\k_\la(v,u_i(v)).$ L'unicit\'e suit ais\'ement.

\qed

\vspace{.1in}  
Nous avons imm\'ediatement le corollaire suivant

\begin{cor}
La dimension de ${\rm{Isom}}\,(\G_\la,\k_\la)$ est donn\'ee par
$$
3n+2+2\sum_{i=1}^s\, r_i^2,
$$
o\`u $2n+2$ d\'esigne la dimension de $G_\lambda.$
\end{cor}

\vspace{.2in}
Pour chaque $1\le j\le s,$ posons $R_j(t):V_j\to V_j$ 
$$ 
\begin{array}{ccc}
R_j(t)e_k&=&\phantom{-}\cos(t)e_k+\sin(t)\check{e}_k\\
R_j(t)\check{e}_k&=&-\sin(t)e_k+\cos(t)\check{e}_k 
\end{array}
$$
pour $e_k, \check{e}_k\in V_j.$ C'est-\`a-dire, si l'on identifie
$V_j$ et ${\mathbb{C}}^{r_j}$ alors $R_j(t)(z_1,\ldots,z_{r_j})=
({\rm{e}}^{it}z_1,\ldots,{\rm{e}}^{it}z_{r_j}).$
Nous avons,

\begin{pr}\label{pr:2.2}
 Dans le cas du groupe oscillateur $(\G_\la,\k_\la)$ les actions du
Th\'eor\`eme \ref{theo:acc} sont donn\'ees respectivement par         
$$
\begin{array}{lccl}
&\sigma\cdot u&=&\left(v_1,R_1(-\dps\frac{t\la_1}{2})u_1R_1(\dps\frac{t\la_1}{2});
\ldots ;v_s,
R_s(-\dps\frac{t\la_s}{2})u_sR_s(\dps\frac{t\la_s}{2})\right),\\
{\mbox{et}}&&&\\ 
&u\cdot\sigma &=&
\left(t,S(u,\sigma),\ldots, \dps\frac{2}{\la_j}\sin(\dps\frac{t\la_j}{2})R_j(\dps\frac{t\la_j}{2})v_j+u_jz_j,
\ldots \right)
\end{array}
$$
o\`u $\sigma=(t,s,z_1,\ldots ,z_n),$ 
$u=(v_1,u_1;\ldots ;v_s,u_s)$ et 
$$S(u,\sigma)=
s-\sum_j\k\left(v_j,\frac{\sin (t\la_j)}{2\la_j}v_j+
\cos (\frac{t\la_j}{2})u_j\circ R_j(\frac{-t\la_j}{2})z_j\right).$$
\end{pr}

\vspace{.1in}
\noindent
{\bf{Preuve.}} Dans le cas des groupes orthorgonaux, l'exponentielle du 
groupe et l'exponentielle de la m\'etrique co\"\i ncident. Or, pour 
$(\G_\la,k_\la),$ cette exponentielle est donn\'ee par la formule:
$$
\exp(t,s,z_i,\ldots,z_n)=$$
$$
(t,s+\frac{1}{2}\sum_j\|z_j\|^2\frac{t\la_j-\sin t\la_j}{t^2\la_j^2},
\frac{{\rm{e}}^{it\la_1}-1}{it\la_1}z_1,\ldots,
\frac{{\rm{e}}^{it\la_n}-1}{it\la_n}z_n),
$$
et, par cons\'equent, la polaire associ\'ee \`a $u=(v_1,u_1;\ldots;v_n,u_n)$ 
est donn\'ee par
$$
{\cal{P}}_u (T,S,\ldots, Z_j,\ldots)=
$$
$$
(T,S-\sum_j k_j(v_j, A_j),\ldots,
\frac{2}{\la_j}\sin\frac{T\la_j}{2}R_j(\frac{T\la_j}{2})v_j+
R_j(\frac{T\la_j}{2})u_jR_j(\frac{-T\la_j}{2})Z_j,\ldots),
$$
o\`u $A_j=(\sin T\la_j/2\la_j)vj+\cos (T\la_j/2)u_jR_j(-T\la_j/2)Z_j.$

\qed

\begin{pr}\label{pr:2.1}
Soit $\la=(\lambda_1,\ldots,\lambda_n)$ tel que $\lambda_1<\cdots<\lambda_n.$ Alors  l'action
de $\G_\la$ sur ${\cal{O}}_{{\cal{R}}}(\g_\la,\k_\la)$ est l'action
triviale, l'action de ${\cal{O}}_{{\cal{R}}}(\g_\la,\k_\la)$ sur $\G_\la$ se fait
par des automorphismes du groupe, i.e. les polaires sont des automorphismes de
$\G_\la$ et le groupe $\left({\rm{Isom}}\,(\G_\la,\k_\la)\right)_0$ est 
isomorphe au groupe produit semi-direct 
$\G_\la\rtimes\left({\cal{O}}_{{\cal{R}}}(\g,\k)\right)_0:$
$$
\left({\rm{Isom}}\,(\G_\la,\k_\la)\right)_0\stackrel{\sim}{=}
\G_\la\rtimes\left({\cal{O}}_{{\cal{R}}}(\g,\k)\right)_0.
$$
\end{pr}

\vspace{.1in}
\noindent
{\bf{Preuve.}} Les sous-espaces $V_i$ \'etant de dimension deux, 
${\cal{O}}(V_i,\k_i)$ est commutatif et, par cons\'equent,
$R_i(-\dps\frac{t\la_i}{2})u_iR_i(\dps\frac{t\la_i}{2})=u_i,$
d'o\`u la premi\`ere assertion. L'\'equation (\ref{eq:pol}) devient
$$
{\cal{P}}_u\circ L_{\sigma}=L_{{\cal{P}}_u(\sigma)}\circ {\cal{P}}_{u}
$$
ce qui revient \`a dire que ${\cal{P}}_u$ est un automorphisme du groupe
$\G_\la.$
\qed

\vspace{.1in}
\noindent
{\bf{Remarque.}}
Sous l'hypoth\`ese de la proposition \ref{pr:2.1}, le groupe
$\left({\rm{Isom}}\,(\G_\la,\k_\la)\right)_0$ s'exprime aussi comme un produit semi-direct
$$
\left({\rm{Isom}}\,(\G_\la,\k_\la)\right)_0\stackrel{\sim}{=}
\left(\G_\la\times V_1\times\cdots\times V_n\right)\rtimes
\left(S^1\right)^n,
$$
o\`u l'action de la droite sur 
$V_1\times\cdots\times V_n\times V_1\times\cdots\times V_n$ est donn\'ee
par
$$
t\cdot(z_1, \ldots,z_n;z_1', \ldots,z_n')
=(\exp(i\la_1t)z_1, \ldots,\exp(i\la_nt)z_n;z_1', \ldots,z_n').
$$
Ce groupe peut-\^etre vu comme un un groupe oscillateur ``d\'eg\'en\'er\'e". 
Cependant il n'est
pas orthogonal.

\begin{pr}\label{pr:laiguales}
Soit $\la=(\lambda_1,\ldots,\lambda_n)$ tel que $\lambda_1=\cdots=\lambda_n,$
$n>2.$ 
Alors  ${\cal{O}}_{{\cal{R}}}(\g_\la,\k_\la)$ est le groupe des 
d\'eplacements
rigides de $(V:=\g'/{\cal{Z}}(\g), \tilde{\k})$ o\`u  $\tilde{\k}$ est la m\'etrique
induite sur $V$ par $\k_\la,$ et $\G_\la,$ identifi\'e au sous-groupe 
$\G_\la\times\{(0,{\rm{Id}})\}$ de
$\left({\rm{Isom}}\,(\G_\la,\k_\la)\right)_0,$ n'est pas distingu\'e
dans $\left({\rm{Isom}}\,(\G_\la,\k_\la)\right)_0.$
\end{pr}

\vspace{.1in}
\noindent
{\bf{Preuve.}}
Les actions du Th\'eor\`eme \ref{theo:acc} \'etant donn\'ees par
$$
\sigma\cdot u=\left(v,R(-\frac{t}{2})\nu R(\frac{t}{2})\right),
$$
et 
$$
u\cdot\sigma=
\left(t,s-\k(v,\frac{\sin (t)}{2}v
+\cos (\frac{t}{2})\nu\circ R(\frac{t}{2})z),2\sin(\frac{t}{2})R(\frac{t}{2})v+\nu (z)
 \right)
$$
o\`u $\sigma=(t,s,z)$ et 
$u=(v,\nu),$ on v\'erifie que la projection de
$
(\sigma_1, u_1)\cdot(\sigma_2, u_2)\cdot(\sigma_1, u_1)^{-1}$ sur
${\cal{O}}_{{\cal{R}}}(\g_\la,\k_\la)$ est \'egale
\`a:
$$
\left(v_1-R(\frac{t}{2})\nu_1R(-t)\nu_1^{-1}R(\frac{t}{2})v_1, 
R(\frac{t}{2})\nu_1R(-t)\nu_1^{-1}R(\frac{t}{2})\right)\ne (0,{\rm{Id}}).
$$
\qed

\section{ G\'eom\'etrie des vari\'et\'es oscillatrices}

    D\'esignons par $\Gamma$ un groupe qui op\`ere \`a droite sur un groupe de 
Lie orthogonal $G$ librement et de fa\c{c}on proprement discontinue et par $p$ 
la projection canonique de $G$ sur $M=\G/\Gamma.$ Si $g'$ est une m\'etrique
pseudo-riemannienne sur $M,$ il existe une unique m\'etrique $g$ sur $G$ telle 
que $p$
est un rev\^etement pseudo-riemannien. La m\'etrique $g$ est invariante par 
$\Gamma.$ Si
de plus l'action naturelle de $G$ sur $M$ se fait par des isom\'etries, $g$ est
invariante par l'action par multiplications \`a gauche de $G$ sur lui m\^eme.
En g\'en\'eral, $g$ n'est pas bi-invariante. On a par exemple

\begin{pr}
La vari\'et\'e oscillatrice de dimension 4, $G_1/\Gamma,$ o\`u 
$$\Gamma=\{(2\pi n, s, a, b); n,2s,a,b\in \mathbb{Z}\}$$ est munie de la 
m\'etrique lorentzienne transitive (donc compl\`ete) de m\'etrique
relev\'ee sur $G_1$ non bi-invariante donn\'ee par la forme quadratique sur 
${\rm{L}}(G_1)$ donn\'ee par
$
<x,x>:= \alpha x_{-1}^2+2x_{-1}x_0+x_1^2+x_2^2.
$
\end{pr}

\vspace{.1in}
Soient $(G_\lambda,\k_\lambda)$ un groupe oscillateur vu comme vari\'et\'e
lorentzienne sy\-m\'e\-tri\-que et $M_\lambda:=G_\lambda/\Gamma $ 
une vari\'et\'e oscillatrice. La
m\'etrique $\k_\lambda$ d\'etermine une m\'etrique sur $M_\lambda$ telle 
que la projection $p:G_\la\to M_\la$ est un
rev\^etement lorentzien. Cette m\'etrique est localement sym\'etrique.
Comme cette derni\`ere est transitive il suit d'un r\'esultat d\^u \`a 
Marsden quelle est compl\`ete. Par cons\'equent elle est sym\'etrique et 
\'evidemment lorentzienne. Ces propri\'et\'es des exemples de \cite{kn:MR1} 
se retrouvent dans des cas plus g\'en\'eraux comme le montre le r\'esultat 
suivant.

\begin{theo}\label{theo:todo}
 Quelle que soit la vari\'et\'e oscillatrice $M_\lambda$ elle est 
 munie:
\begin{enumerate}
\item D'une structure de vari\'et\'e affine dont la preimage par $p$ 
est la structure affine sur $G_\lambda$ du th\'eor\`eme 4.2,
\item De m\'etriques pseudo-riemanniennes localement sym\'etriques compl\`etes,
en g\'en\'eral non transitives, induites par la famille des m\'etriques sur 
$G_\lambda$ d\'e\-cri\-tes dans le th\'eor\`eme 4.3,
\item De m\'etriques non compl\`etes ( et donc non transitives) d'indice 
quelconque non nul telles que $p$ est pseudo-riemannienne relativement aux 
m\'etriques sur $G_\lambda$ de la proposition \ref{pr:ind}.
\end{enumerate}
\end{theo}

\begin{pr}\label{pr:ind}
Les groupes oscillateurs admettent des m\'e\-tri\-ques d'indice quelconque 
non nul  non compl\`etes et donc non transitives.
\end{pr}

\vspace{.1in}
\noindent
{\bf{Preuve.}}
Consid\'erons d'abord sur le groupe oscillateur de dimension 4
les m\'etriques associ\'ees aux isomorphismes lin\'eaires
$\k$-sym\'etriques suivants
$$
\begin{array}{rclrclrclrcl}
u_1(e_{-1})&=&e_1&u_1(e_{0})&=&e_{-1}&u_1(e_{1})&=&e_0&u_1(\check{e}_{1})&=&\check{e}_1\\
u_2(e_{-1})&=&\check{e}_1&u_2(e_{0})&=&e_1&u_2(e_{1})&=&e_{-1}&u_2(\check{e}_{1})&=&e_0.
\end{array}
$$

La premi\`ere m\'etrique est d'indice 1, la seconde
d'indice 2. Montrons qu'elles dont incompl\`etes, c'est-\`a-dire
que les champs de vecteurs \ref{eq:geo2} n'est pas complet.
En ce qui concerne $u_1$ on constate que 
$$
\begin{array}{lcl}
\gamma_1(t)&=&(c,c,c-(2\rho^2/c)\sec^2(\rho\, t),
-2\rho\tan(\rho\, t))
\end{array}
$$
(o\`u $\rho$ est un r\'e\'el non nul quelconque)
est une courbe
int\'egrale de (\ref{eq:geo2}) qui n'est pas compl\`ete.
Dans le cas de $u_2,$ l'\'equation (\ref{eq:geo2})
s'\'ecrit
$$
\begin{array}{cclcccl}
\dot{x}_{-1}&=&-x_{-1}x_0+x_1^2&&\dot{x}_{0}&=&-x_{1}\check{x}_1+x_{-1}^2\\
\dot{x}_{1}&=&\phantom{-}0&&\dot{\check{x}}_{1}&=&-x_1x_{-1}+x_{0}\check{x}_1.
\end{array}
$$
Celle-ci admet comme int\'egrales premi\`eres 
$2x_1\check{x}_1+x_{-1}^2+x_{0}^2$ et 
$ x_{-1}\check{x}_1+x_{0}x_1.$
L'incompl\'etude de ce champ est cons\'equence de la non-compl\'etude, dans la
direction d'un vecteur de type temps, de l'\'equation
$\dot{x}=x^2/2,$ v\'erifiant la condition initiale 
$x(0)=x_0.$ 

Soit $v_i$ un isomorphisme $\k_i$-sym\'etrique de 
${\rm{Vect}}\,\{e_i,\check{e}_i\}$ pour tout $i=2,\ldots,n.$ Alors 
$u_i\oplus v_2\oplus \cdots \oplus v_n$ est un isomorphisme $\k_\la$-sym\'etrique 
de $\g_\la.$

Si $\gamma_i$ est une trajectoire de (\ref{eq:geo2}) relativement \`a $u_i,$ 
il est clair que
$(\gamma_i,0,\ldots ,0)$ une trajectoire de (\ref{eq:geo2}) dans $\g_\la,$
relativement \`a 
$u_i\oplus v_2\oplus \cdots \oplus v_n.$ \'Evidemment les $v_j$ 
peuvent \^etre choisis de sorte que la m\'etrique qui leur est associ\'ee soit d'indice 
arbitraire non nul .
\qed

\end{document}